\documentclass{article}
\usepackage{amsthm}

\usepackage[T1]{fontenc} %accents
\usepackage{amsmath}
\usepackage{amscd}
\usepackage{amssymb}
\usepackage{latexsym}
\usepackage{xypic}

\newcommand{\LL}{{\cal L}}

\newcommand{\CA}{{\cal A}}

\newcommand{\CK}{{\cal K}}

\newcommand{\FA}{F(\mathcal{A})}

\newcommand{\chop}{\dag}

\def\epsilon{\varepsilon}

\newcommand{\Out}{\mbox{Out}}

\newcommand{\inv}{^{-1}}

\newcommand{\dom}{\mbox{dom}}
\newcommand{\Adm}{\mbox{Adm}}

\newcommand{\FN}{F_N}   % F ou F_n ou F_N ?
\newcommand{\cvn}{\mbox{cv}_N}   
\newcommand{\CVN}{\mbox{CV}_N}   
\newcommand{\barCVN}{\bar{\mbox{CV}}_N}   
\newcommand{\barcvn}{\bar{\mbox{cv}}_N}   
\newcommand{\CQ}{{\cal Q}}
\newcommand{\Tobs}{\widehat T^{\mbox{\scriptsize obs}}}

\newcommand{\CLadm}{\mathcal{L}_{\mbox{\scriptsize adm}}}

\newcommand{\Loneadm}{L^1_{\mbox{\scriptsize adm}}}   
\newcommand{\Ladm}{L_{\mbox{\scriptsize adm}}}   
\newcommand{\CQT}{\CQ_T}
\newcommand{\CQK}{\CQ_\CK}
\newcommand{\TKmin}{T_\CK^{\mbox{\scriptsize{min}}}}
\newcommand{\jmin}{j^{\mbox{\scriptsize{min}}}}

\newcommand{\R}{\mathbb R}
\newcommand{\Z}{\mathbb Z}

\newcommand{\N}{\mathbb N}

\renewcommand{\L}{{\cal L}}

\def\strutdepth{\dp\strutbox}
\def \ss{\strut\vadjust{\kern-\strutdepth \sss}}
\def \sss{\vtop to \strutdepth{
\baselineskip\strutdepth\vss\llap{$\diamondsuit\;\;$}\null}}

\def\strutdepth{\dp\strutbox}
\def \sst{\strut\vadjust{\kern-\strutdepth \ssss}}
\def \ssss{\vtop to \strutdepth{
\baselineskip\strutdepth\vss\llap{$\spadesuit\;\;$}\null}}

\def\strutdepth{\dp\strutbox}
\def \ssh{\strut\vadjust{\kern-\strutdepth \sssh}}
\def \sssh{\vtop to \strutdepth{
\baselineskip\strutdepth\vss\llap{$\heartsuit\;\;$}\null}}

\def\qed{\hfill\rlap{$\sqcup$}$\sqcap$\par}
\def\bar{\overline}
\def\tilde{\widetilde}
\def\hat{\widehat}

\newtheorem{thm}{Theorem}[section]
\newtheorem{cor}[thm]{Corollary}
\newtheorem{lem}[thm]{Lemma}
\newtheorem{prop}[thm]{Proposition}

\newtheorem*{cor61*}{Corollary~\ref{cor:geometrictree}}
\newtheorem*{cor63*}{Corollary~\ref{cor:finite-bounded-approximation}}
\theoremstyle{definition}
\newtheorem{defn}[thm]{Definition}

\newtheorem{rem}[thm]{Remark}

\begin{document}

\title{$\R$-trees, dual laminations, and compact systems of partial isometries%\\
%Arbres r\'eels, laminations duales et syst\`emes compacts d'isom\'etries partielles
}

% \author[Thierry Coulbois, Arnaud Hilion and Martin Lustig]{Thierry Coulbois, Arnaud Hilion and Martin Lustig \\
% LATP, universit\'e Aix-Marseille III, Marseille, France
% }

\author{Thierry Coulbois, Arnaud Hilion and Martin Lustig \\
LATP, universit\'e Aix-Marseille III, Marseille, France
}

\date{\today }

\maketitle

\begin{abstract}
  Let $\FN$ be a free group of finite rank $N \geq 2$, and let $T$ be
  an $\R$-tree with a very small, minimal action of $\FN$ with dense
  orbits. For any basis $\CA$ of $\FN$ there exists a {\em heart}
  $K_{\CA} \subset \bar T$ (= the metric completion of $T$) which is a
  compact subtree that has the property that the dynamical system of
  partial isometries $a_{i} : K_{\CA} \cap a_{i} K_{\CA} \to a_{i}\inv
  K_{\CA} \cap K_{\CA}$, for each $a_{i} \in \CA$, defines a tree
  $T_{(K_{\CA}, \CA)}$ which contains an isometric copy of $T$ as
  minimal subtree.
\end{abstract}

%% \begin{abstract}
%% Soit $\FN$ un groupe libre de rang fini $N\geq 2$, et soit $T$ un
%% arbre r\'eel munit d'une action minimale, tr\`es petite et \`a orbites
%% denses de $\FN$ par isom\'etries. Pour chaque base $\CA$ de $\FN$ il
%% existe un {\em c\oe{}ur} $K_{\CA} \subset \bar T$ (le compl\'et\'e
%% m\'etrique de $T$) qui est un arbre compact tel que le syt\`eme
%% dynamique form\'e des ism\'etries partielles $a_{i} : K_{\CA} \cap
%% a_{i} K_{\CA} \to a_{i}\inv K_{\CA} \cap K_{\CA}$, pour chaque
%% $a_i\in\CA$ , definit un arbre $T_{(K_{\CA}, \CA)}$ qui contient une
%% copie isom\'etrique de $T$ comme sous-arbre minimal.
%% \end{abstract}

%\tableofcontents

\section{Introduction}
\label{introduction}

A point on Thurston's boundary of Teichm\"uller space ${\cal T}(\Sigma)$
for a surface $\Sigma$ can be understood alternatively as a measured geodesic
lamination $({\mathfrak L}, \mu)$ on $\Sigma$, up to rescaling of the
transverse measure, or as a small action of $\pi_{1} \Sigma$ on some
$\R$-tree $T$, up to $\pi_{1} \Sigma$-equivariant homothety. The
correspondence between these two objects, which are naturally dual to
each other, is given by the fact that points of $T$ are in 1-1
correspondence (or ``one-to-finite'' correspondence, for the
branchpoints of $T$) with the leaves of $\tilde {\mathfrak L}$, i.e. the
lift of $\mathfrak L$ to the universal covering $\tilde \Sigma$.  The
metric on $T$ is determined by $\mu$, and vice versa.

\smallskip

Culler-Vogtmann's Outer space $\CVN$ is the analogue of ${\cal
T}(\Sigma)$, with $\Out(\FN)$ replacing the mapping class group.  A
point of the boundary $\partial \CVN$ is given by a homothety class
$[T]$ of very small isometric actions of the free group $\FN$ on an
$\R$-tree $T$.  In general, $T$ will not be dual to a measured
lamination on a surface.  However, in \cite{chl1-I, chl1-II} an
``abstract'' {\em dual lamination} $L(T)$ has been defined for any
such $T$, which is very much the analogue of $\mathfrak L$ in the
surface case. The dual lamination $L(T)$ is an algebraic lamination:
it lives in the double Gromov boundary of $\FN$, and the choice of a
basis $\CA$ transforms $L(T)$ into a symbolic dynamical system which
is a classical subshift in $\CA\cup\CA^{-1}$.  The dual lamination
$L(T)$, and variations of it, have already been proved to be a useful
invariant of the tree $T$, compare \cite{bfh, chl2, hm, kl2}.

In the case of measured laminations on a surface, the standard tool
which allows a transition from geometry to combinatorial dynamics, is
given by interval exchange transformations.  The combinatorics which
occur here are classically given through coding geodesics on a surface
by sequences of symbols, where the symbols correspond to subintervals,
and the sequences are given by the first return map.  Conversely, the
surface and the lamination (or rather ``foliation'', in this case),
can be recovered from the interval exchange transformation by suspension,
i.e. by realizing the map which exchanges the subintervals by a
(foliated) mapping torus.

\smallskip

Taking the basic concept of this classical method one step further and
considering directly the dual tree $T$ rather than the lamination
given by the combinatorial data, one considers for any $[T] \in
\partial \CVN$ a finite metric subtree $K \subset T$, and for some
basis $\CA$ of $\FN$ the induced finite system of partial isometries
between subtrees of $K$: Each basis element $a_{i} \in \CA$ defines a
partial isomerty $a_{i} : K \cap a_{i}K \to a_{i}\inv K \cap K$, and
these partial isometries play the role of the interval exchange
transformation.  Any such pair $\CK = (K, \CA)$ gives canonically rise
to a tree $T_{\CK}$ together with an $\FN$-equivariant map $j: T_{\CK}
\to T$. The tree $T_{\CK}$ is the ``unfolding space'' of the system
$\CK$.  The class of $\R$-trees $T$, with the property that for some
such finite $K$ the map $j$ is an isometry, have been investigated
intensely, and they play an important role in the study of $\partial
\CVN$, see \cite{gl}.

\smallskip

Indeed, if $K$ is an interval and if it is simultaneously equal to the
union of domains and the union of ranges of the isometries (and if
these unions are disjoint unions except at the boundary points), then
$\CK$ defines actually an interval exchange transformation.  If one only
assumes that $K$ is finite, this will in general not be true: one only
obtains a system of interval translations (see for instance
\cite{bh}).  On the level of $\R$-trees one obtains in the first case
{\em surface tree actions}, and in the second case actions that were
alternatively termed {\em Levitt, thin} or {\em exotic}.  The union of
these two classes are precisely the actions called {\em geometric} in
\cite{gl}.

\smallskip

However, both of these types of actions seem to be more the exception
than the rule: Given any point $[T] \in \partial \CVN$, there is in
general no reason why $T$ should be determined by a system of partial
isometries based on a finite tree $K \subset T$. A possible way to
deal with such $T$ is to consider increasing sequences of finite
subtrees and thus to approximate $T$ by the sequence of ensuing
geometric trees $T_{\CK}$, in the spirit of the ``Rips machine'',
which is an important tool to analyze arbitrary group actions on
$\R$-trees.  The goal of this paper is to propose a more direct
alternative to this approximation technology:

\smallskip

We replace the condition on the subtree $K \subset T$ to be finite by
the weaker condition that $K$ be compact.  It turns out that almost
all of the classical machinery developed for the approximation trees
$T_{\CK}$ for finite $K$ carries over directly to the case of compact
$K$.  However, the applications of such $T_{\CK}$ concern a much larger
class of trees: In particular, every minimal very small $T$ with dense
orbits can be described directly, i.e. circumventing completely the
above approximation, as minimal subtree $T^{min}_{\CK}$ of the tree
$T_{\CK}$, for a properly chosen compact subtree $K$ of the metric
completion $\bar T$ of $T$.

\begin{thm}\label{thm:main}
  Let $T$ be an $\R$-tree provided with a very small, minimal,
  isometric action of the free group $\FN$ with dense orbits.  Let
  $\CA$ be a basis of $\FN$.  Then there exists a unique compact subtree
  $K_\CA \subset \bar T$ (called the ``heart'' of $T$ w.r.t. $\CA$),
  such that for any compact subtree $K$ of $\bar T$ one has:
\[
T=T^{min}_\CK  \iff  K_\CA\subseteq K
\]
\end{thm}

This is a slightly simplified version of Theorem \ref{thm:mainthm}
proved in this paper.  The main tool for this proof (and indeed for
the definition of the heart $K_{\CA}$) is the dual lamination $L(T)$.
We define in this article (see \S\ref{sectionthree}) a second {\em
admissible lamination} $\Ladm(\CK)$ associated to the system of
partial isometries $\CK = (K, \CA)$.  One key ingredient in the
equivalence of Theorem \ref{thm:main} is to prove that the two
statements given there are equivalent to the equation $L(T) =
\Ladm(\CK)$.  The other key ingredient, developed in
\S\ref{sectionfour}, is a new understanding of the crucial map $\CQ:
\partial \FN \to \bar T \cup \partial T$ from \cite{ll4}, based on the
dynamical system $\CK = (K, \CA)$.  The proof of Theorem
\ref{thm:mainthm} uses the full strength of the duality between trees
and laminations, and in particular a transition between the two, given
by the main result of our earlier paper \cite{chl2}.

\medskip

We would like to emphasize that the main object of this paper, the
{heart} $K_{\CA}$ of $T$ with respect to any basis $\CA$ of $\FN$, is
a compact subtree of $\bar T$ that is determined by algebraic data
associated to $T$, namely by the dual algebraic lamination $L(T)$ of
$T$.  This system $\CK_{\CA} = (K_{\CA}, \CA)$ of partial isometries
is entirely determined by the choice of the basis $\CA$ and it depends
on $\CA$, but important properties of it turn out to be independent of
that choice.  For example, we derive in \S\ref{sec:applications} from
the above theorem the following direct characterization of geometric
trees, and we also give a sharpening of Gaboriau-Levitt's
approximation result for trees $T$ from $\partial \CVN$:

\begin{cor61*}
A very small minimal $\R$-tree $T$, with isometric $\FN$-action that
has dense orbits, is geometric if and only if, for any basis $\CA$ of
$\FN$, the heart $K_{\CA} $ is a finite subtree of $T$.
\end{cor61*}

\begin{cor63*}
For every very small minimal $\R$-tree $T$, with isometric
$\FN$-action that has dense orbits, there exists a sequence of finite
subtrees $K(n)$ of uniformely bounded diameter, such that:
\[
T = \lim_{n \to \infty} \, T_{\CK(n)}
\]
\end{cor63*}

In contrast to the case of geometric $\R$-trees, there are trees in
$\partial \CVN$ for which the compact heart is far from being
finite. Indeed it is proven in \cite{coulbois-fractal} that the
compact heart of the repulsive tree $T_{\Phi\inv}$ of an iwip outer
automorphism $\Phi$ of $\FN$ has Hausdorff dimension equal to
$\max(1,\frac{\ln\lambda_\Phi}{\ln\lambda_{\Phi\inv}})$, where
$\lambda_\Phi$ and $\lambda_{\Phi\inv}$ are the expansion factors of
$\Phi$ and $\Phi\inv$ respectively. As these expansion factors are in
general not equal, we can assume that
$\lambda_\Phi>\lambda_{\Phi\inv}$ to get a compact heart with
Hausdorff dimension strictly bigger than $1$.

\section{$F_{N}$-actions on $\R$-trees and their heart}
\label{sectiontwo}

In this section we first recall some well known facts about $\R$-trees
$T$ with isometric action of a free group $\FN$. We also recall
algebraic laminations, and in particular the dual lamination $L(T)$.
We then concentrate on the specific case of very small trees with
dense orbits, and for such trees we define the limit set and the heart
of $T$ with respect to a fixed basis $\CA$ of $\FN$.

In this paper we need some of the machinary developed in our previous
articles \cite{chl1-I, chl1-II, chl2}. We present these tools in this
section, but refer to those articles for proofs and for a more
complete discussion.

\subsection{Background on $\R$-trees} 
\label{prelims}

An $\R$-tree $T$ is a metric space which is $0$-hyperbolic and
geodesic.  Alternatively, a metric space $T$ is an $\R$-tree if and
only if any two points $x,y \in T$ are joined by a unique topological arc $[x,y]
\subset T$, and this arc (called a {\em segment}) is geodesic.  For
any $\R$-tree $T$, we denote by $\bar T$ the metric completion and by
$\partial T$ the Gromov boundary of $T$. We also write $\hat T = \bar
T\cup\partial T$.

\smallskip

Most $\R$-trees $T$ considered in this paper are provided with an
action by isometries (from the left) of a non-abelian free group $\FN$
of finite rank $N\geq 2$. Such an action is called {\em minimal} if $T$
agrees with its minimal $\FN$-invariant subtree. We say that the
action has {\em dense orbits} if for some (and hence every) point $x \in
T$ the orbit $\FN\cdot x$ is dense in $T$.  In the case of dense orbits,
the following three conditions are equivalent:
\begin{itemize}
\item $T$ has {\em trivial arc stabilizers} (i.e. for any distinct $x,
  y \in T$ and $w \in \FN$ the equality $w[x,y] = [x, y]$ implies $w =
  1$).

\item The $\FN$-action on $T$ is {\em small} (see \cite{cm,chl1-II}).

\item The $\FN$-action on $T$ is {\em very small} (see
  \cite{cl,chl1-II}).

\end{itemize}

As usual, for any $w \in \FN$ we denote by $\|w\|_{T}$ (or simply
by $\|w\|$) the translation length of the action of $w$ on $T$, i.e. the
infimum of $d(x, wx)$ over all $x \in T$.

There are two types of isometries of $T$:  An element $w \in 
\FN$ acts as an {\em elliptic} isometry on $T$ if it fixes a point, which is equivalent to  
$\|w\| = 0$.  If $\|w\| > 0$, then the action of $w$ on $T$ is 
called {\em hyperbolic}:  There is a well defined {\em axis} in $T$, 
which is isometric to $\R$ and is $w$-invariant: the element $w$ translates 
every point on the axis by $\|w\|$.

\smallskip

A continuous map $T \to T'$ between $\R$-trees is called a {\em
  morphism} if every segment is mapped locally injectively except at
finitely many points.

\subsection{The observers' topology on $T$}
\label{observerstopology}

There are various independent approaches in the literature to define
$\R$-trees as topological spaces without reference to the metric.  The
following version has been studied in \cite{chl2}.

\begin{defn}\label{def:obs} 
Let $T$ be an $\R$-tree.  A {\em direction} in $\hat T$ is a
connected component of the complement of a point of $\hat T$.  A subbasis
of open sets for the {\em observers' topology} on $\hat T$ is given by the set 
of all such
directions in $\hat T$.
\end{defn}

The observers' topology on $\hat T$ (or $T$) is weaker than the metric
topology: For example, any sequence of points that ``turns around'' a
branch point converges to this branch point.  We denote by $\Tobs$ the
set $\hat T$ equipped with the observers' topology. The space $\Tobs$
is Hausdorff and compact.

\medskip

For any sequence of points $P_{n}$ in $\hat T$, and for some
{\em base point} $Q \in \hat T$,
there is a well defined {\em inferior
limit from $Q$},
which we denote
by: 
$$P = \underset{n \to \infty}{\liminf}{}_{Q} \,\, P_{n}$$
It is given by
$$[Q, P]=\overline{\bigcup_{m=0}^\infty\bigcap_{n\geq m}[Q, P_n]}.$$

The inferior limit $P$ is always contained in the closure of the
convex hull of the $P_{n}$, but its precise location does in fact
depend on the choice of the base point $Q$.  However, in \cite{chl2}
the following has been shown:

\begin{lem}
\label{lem:convergence}
If a sequence of points $P_{n}$ converges in $\Tobs$ to some limit
point $P \in \Tobs$, then for any $Q \in \hat T$ one has:
\[
P = \underset{n \to \infty}{\liminf}{}_{Q} \,\, P_{n}
\]
\end{lem}

The observers' topology is very useful, but it is also easy to be 
deceived by it.  For example, it is not true that any continuous map 
between $\R$-trees $T_{1} \to T_{2}$ induces canonically a continuous map 
$\Tobs_{1} \to \Tobs_{2}$, as is illustrated in the following remark.

\begin{rem}
\label{notcontinuous}
Let $T_{1}$ be the $\infty$-pod, given by a center $Q$ and edges 
$[Q, P_{k}]$ of length 1, for every $k \in \N$. Let $T_{2}$ be 
obtained from $T_{1}$ by gluing the initial segment of length 
$\frac{k-1}{k}$ of each $[Q, P_{k}]$, for $k \geq 2$, to $[Q, P_{1}]$.
Then the canonical map $f: T_{1} \to T_{2}$ is continuous, and even a 
length decreasing morphism, but $\lim P_{k} = Q$, while $\lim 
f(P_{k}) = f(P_{1}) \neq f(Q)$.
\end{rem}

We refer the reader to \cite{chl2} for more details about the observers' 
topology.

\subsection{Algebraic laminations}
\label{algebraiclaminations}

For the free group $\FN$ of finite rank $N \geq 2$, we denote by 
$\partial \FN$ the Gromov boundary of $\FN$. We also consider
$$
\partial^{2}\FN = \partial \FN \times \partial \FN \smallsetminus 
\Delta \, ,
$$
where $\Delta$ denotes the diagonal.  The space $\partial^{2}\FN$ 
inherits from $\partial \FN$ a left-action of $\FN$, defined by $w(X, 
Y) = (wX, wY)$ and a topology. It also admits the {\em flip map} $(X, 
Y) \mapsto (Y, X)$.  An {\em algebraic lamination} $L^{2} \subset 
\partial^{2}\FN$ is a non-empty closed subset which is invariant under the 
$\FN$-action and the flip map. 

\smallskip

If one choses a basis $\CA$ of $\FN$, then every element $w \in \FN$ 
can be uniquely written as a finite reduced word in $\CA^{\pm 1}$, so that 
$\FN$ is canonically identified with the set $\FA$ of such words. 
Similarly, a point of the boundary $\partial \FN$ can be written as an
infinite reduced word $X = z_{1} z_{2} \ldots$, so that $\partial \FN$ 
is canonically identified with the set $\partial \FA$ of such infinite 
words. 

\smallskip

We also consider reduced biinfinite indexed words 
$$Z = \ldots z_{-1} z_{0} z_{1} \ldots$$
with all $z_{i} \in \CA^{\pm 1}$. We say that $Z$
has {\em positive half}
$Z^+ = z_{1} z_{2} \ldots$ and {\em negative half} $Z^- = z_{0}\inv 
z_{-1}\inv \ldots$, which are two infinite words
$$
Z^{+}, Z^- \in \partial \FA
$$
with distinct initial letters $Z^{+}_{1} \neq Z^-_{1}$.  We write the reduced product $Z =
(Z^-)\inv \cdot Z^{+}$ to mark the letter $Z_1^+$
with index 1.

For any fixed choice of a basis $\CA$, an algebraic lamination $L^{2}$ 
determines a {\em symbolic lamination} 
$$L_{\CA} = \{(Z^-)\inv \cdot Z^{+} \mid (wZ^-, wZ^{+}) \in L^{2}\}$$
as well as a {\em laminary language}
$$\LL_{\CA} = \{w \in \FA \mid w \mbox{ is  a subword 
of  some }  Z \in L_{\CA} \} \, .$$
Both, symbolic laminations and laminary languages can be characterized 
independently, and the natural transition from one to the other and 
back to an algebraic lamination has been established with care in 
\cite{chl1-I}.  In case we do not want to specify which of the three 
equivalent terminologies is meant, we simply speak of a {\em 
lamination} and denote it by $L$.

\medskip

One of the crucial points of the encounter between symbolic dynamics
and geometric group theory, in the subject treated in this paper,
occurs precisely at the transition between algebraic and symbolic
laminations.  Since the main thrust of this paper (as presented in
\S\ref{sectionthree}) can be reinterpreted as translating the
symbolic dynamics viewpoint into the world of $\R$-trees, it seems
useful to highlight this transition in the symbolic language, before
embroiling it with the topology of $\R$-trees:

\begin{rem}
\label{symbolicpartialisom}
As before, we fix a basis $\cal A$ of $\FN$, and denote an element $X$
of the boundary $\partial\FN = \partial F(\CA)$
by the corresponding infinite reduced word in $\CA^{\pm 1}$. We
denote by $X_n$ its prefix of length $n$.

We consider the {\em unit cylinder} $C^{2}_\CA$ in $\partial^2\FN \,$:
\[
C^{2}_\CA=\{(X,Y)\in\partial^2\FN\ |\ X_1\neq Y_1\}
\]
Contrary to  $\partial^2\FN$, the unit cylinder $C^{2}_\CA$
is a compact set (in fact, a Cantor set).  The unit cylinder
$C^{2}_\CA$ has the property that the canonical map $\rho_{\CA}: (X,
Y) \mapsto X^{-1} \cdot Y$ (see \cite[Remark~4.3]{chl1-I}) restricts to an injection
on $C^{2}_\CA$  with inverse map $Z \mapsto (Z^-,
Z^+)$.

\smallskip

In symbolic dynamics, the natural operator on biinfinite sequences is
the {\em shift} map, which in our notation is given by
$$\sigma(X^{-1} \cdot Y) = X^{-1} Y_{1} \cdot (Y_{1}^{-1} Y) \, ,$$
i.e. the same symbol sequence as in
$X^{-1} \cdot Y$, but with $Y_{1}$ as letter of index 0.

\smallskip

On the other hand, there is a system of ``partial bijections'' on
$C^{2}_\CA$, given for each $a_{i} \in \CA$ by:
$$a_{i}: C^{2}_\CA \cap {a_{i}}\inv C^2_\CA \to a_{i} C^{2}_\CA \cap C^{2}_\CA$$
A particular feature of this system is that it
``commutes'' via the
map $\rho_{\CA}$ with the shift map $\sigma$ on the set of
biinfinite reduced words:  More precisely, for all $(X, Y) \in
C^{2}_{\CA}$ one has:
$$\rho_{\CA}(Y_{1}^{-1} (X, Y)) =
\sigma(\rho_{\CA}(X, Y) )$$ This transition from group action to the
shift (or more precisely, the converse direction), will be explored in
\S\ref{sectionthree} in detail, with the additional feature that the
topology of compact trees is added on, in the analogous way as
interval exchange transformations are a classical tool to interpret
certain symbolic dynamical systems topologically.
\end{rem}

\subsection{The dual lamination $L(T)$}
\label{duallamination}

In \cite{chl1-II} a {\em dual lamination} $L(T)$ for any isometric
action of a free group $\FN$ on an $\R$-tree $T$ has been introduced
and investigated.  If $T$ is very small and has dense orbits, three
different definitions of $L(T)$ have been given in \cite{chl1-II} and
shown there to be equivalent.  However, as in this paper we can not
always assume that $T$ has dense orbits, it is most convenient to fix
a basis $\CA$ of $\FN$ and to give the general definition of $L(T)$
via its {\em dual laminary language} $\LL_{\CA}(T)$ (see
Definition~4.1 and Remark~4.2 of \cite{chl1-I}), which determines
$L(T)$ and vice versa:
\[
\begin{array}[t]{lll}
\LL_{\CA}(T)& = &\{v\in\FA\ |\ \forall \, \epsilon>0 \, \,  \exists \,
u,w\in\FA~:
\|u \cdot v \cdot w\|_{T}<\epsilon, \\
&& \qquad\qquad \quad\, \, \,   u \cdot v \cdot w \, \, \, 
\mbox{reduced and cyclically reduced} \}
\end{array}
\]

\begin{rem}
\label{dualgivenbyminimal}
It follows directly from this definition that $L(T) = L(T^{min})$, 
where $T^{min}$ denotes the minimal $\FN$-invariant subtree of $T$.
\end{rem}

\subsection{The map $\CQ$}\label{sec:Q}

\begin{thm}[\cite{ll4, ll3}]\label{thm:Qexists}
  Let $T$ be an $\R$-tree with a very small action of $\FN$ by
  isometries that has dense orbits.  Then there exists a surjective
  $\FN$-equivariant map $\CQ:\partial\FN\to\hat T$ which has the
  following property:

  For any sequence of elements $u_n$ of $\FN$ which converges to
  $X\in\partial\FN$ and for any point $P\in T$, if the sequence of
  points $u_nP \in T$ converges (metrically) in $\hat T$ to a point $Q$, then
  $\CQ(X)=Q$.
\end{thm}

Using the properties of a metric topology, we get the following lemma.

\begin{lem}\label{lem:QTandK}
  Let $T$ be an $\R$-tree with a very small action of $\FN$ by
  isometries that has dense orbits.  Let $K$ be a compact (with
  respect to the metric topology) subtree of $\bar T$.  Let $Q$ be a
  point in $K$ and $w_n$ a sequence of elements in $\FN$ which
  converge in $\FN\cup\partial\FN$ to some $X\in\partial\FN$. If for
  all $n$ one has $w_n\inv Q\in K$, then $\CQ(X)=Q$.
\end{lem}
\begin{proof}
  As $K$ is compact, up to passing to a subsequence, we can assume
  that $w_n\inv Q$ converges to a point $P$ in $K$, that is to say
  $\lim_{n\to\infty} d(w_n\inv Q,P)=0$. As the action is isometric,
  we get that $\lim_{n\to\infty} d(Q,w_nP)=0$, i.e. the  $w_nP$
  converge to $Q$. Hence Theorem~\ref{thm:Qexists} gives the desired
  conclusion $\CQ(X)=Q$.
%${}^{}$
\end{proof}

It is crucial for the arguments presented in this paper to remember
that the map $\CQ$ is not continuous with respect to the {\em metric
  topology} on $\hat T$, i.e. the topology given by the metric on $T$.  In
fact, this has been the reason why in \cite{chl2} the weaker
observers' topology on $\hat T$ has been investigated.

\begin{thm}[{\cite[Remark 2.2 and Proposition 2.3]{chl2}}] 
\label{Qunique}
Let $T$ be an $\R$-tree with isometric very small action of $\FN$
that has dense orbits. Then the following holds:
\begin{enumerate}
\item[(1)] The map $\CQ$ defined in Theorem \ref{thm:Qexists} is
  continuous with respect to the observers' topology, i.e. it defines
  a continuous equivariant surjection
$$\CQ:\partial\FN\to\Tobs.$$
\item[(2)] For any point $P \in T$ the map $\CQ$ defines the
  continuous extension to $\FN \cup \partial \FN$ of the map
$$\CQ_{P}: \FN \to \Tobs,  \, \, w \mapsto wP\, .$$
\end{enumerate}
\end{thm}

Though obvious it is worth noting that the last property determines the map $\CQ$ 
uniquely.

\subsection{The map $\CQ^2$}
\label{themapQ2}

If the tree $T$ is very small and has dense orbits, the dual
lamination $L(T)$ described in \S\ref{duallamination} admits an
alternative second definition via the above defined map $\CQ$ as
algebraic lamination $L^{2}(T)$ (compare \S\ref{algebraiclaminations}):
\[
L^{2}(T)=\{(X,Y)\in\partial^2\FN\ |\ \CQ(X)=\CQ(Y)\}
\]
It has been proved in \cite{ll4, ll3} that the map $\CQ$ is one-to-one
on the preimage of the Gromov boundary $\partial T$ of $T$.  Hence the
map $\CQ$ induces a map $\CQ^2$ from $L^{2}(T)$ to $\bar T$, given by:
$$\CQ^{2}((X, Y)) = \CQ(X) =\CQ(Y)$$
In light of the above discussion the following result seems 
remarkable. It is also crucial for the definition of the heart of 
$T$ in the next subsection.

\begin{prop}[{\cite[Proposition 8.3]{chl1-II}}]
\label{continuous}
The  
$\FN$-equivariant 
map
$$\CQ^2:L^{2}(T)\to\bar T$$
is continuous, with respect to the 
metric topology on $\bar T$.
\end{prop}

%%%%%%%%%%%%%%%%%%%%
As in \cite[\S 2]{chl2}, we consider
the equivalence relation on $\partial F_N$ whose classes are
fibers of $\CQ$, and we denote by $\partial F_N/L^2(T)$ the quotient
set.  The quotient topology on $\partial F_N/L^2(T)$ is the finest
topology such that the natural projection $\pi:\partial F_N\to\partial
F_N/L^2(T)$ is continuous.  The map $\CQ$ splits over $\pi$, thus
inducing a map $\varphi: \partial F_N/L^2(T) \to \Tobs$ with $\CQ =
\varphi \circ \pi$.

\begin{thm}[{\cite[Corollary 2.6]{chl2}}]  
\label{thm:homeo}
The map 
$$\varphi: \partial F_N/L^2(T) \to \Tobs$$
is a homeomorphism.
\end{thm}

\subsection{The limit set and the heart of $T$}
\label{limitset}
We consider again the {\em unit
cylinder} $C^{2}_\CA =\{(X,Y)\in\partial^2\FN\ |\ X_1\neq Y_1\}$
in $\partial^2\FN \,$ as defined in Remark 
\ref{symbolicpartialisom}.
The following definition is the crucial innovative tool of this 
paper:

\begin{defn}
\label{limitset-heart}
The {\em limit set} of $T$ with respect to the basis $\CA$ is the set
\[
\Omega_\CA=\CQ^2(C^{2}_\CA\cap L^{2}(T)) \subset \bar T\, .
\] 
The {\em heart} $K_\CA$ of $T$ with respect to
the basis $\CA$ is the convex hull in $\bar T$ of 
the limit set $\Omega_\CA$.
\end{defn}

It is not hard to see that in any $\R$-tree the convex hull of a compact 
set is again compact. Thus we obtain
from 
%the previous 
Proposition \ref{continuous} and Definition \ref{limitset-heart}:
%we 
%get:

\begin{cor}
The limit set $\Omega_\CA$ is a compact subset of $\bar T$. The
heart $K_\CA\subset \bar T$ is a compact $\R$-tree. 
\end{cor}

Note that, while $L^{2}(T)$ does not depend on the choice of the basis
$\CA$, the unit cylinder $C^{2}_\CA$ and thus the limit set and the heart
of $T$ do crucially depend on the choice of $\CA$.

\section{Systems of isometries on compact $\R$-trees}
\label{sectionthree}

In this section we review the basic construction that associates an
$\R$-tree to a system of isometries. This goes back to the seminal
papers of D.~Gaboriau, G.~Levitt,
 and F.~Paulin \cite{glp} and M.~Bestvina and M.~Feighn \cite{bf}, and before them to the
study of surface trees and the work of J.~Morgan and P.~Shalen
\cite{ms}, R.~Skora \cite{skora}, and of course to the 
fundamental work of E.~Rips.

\smallskip

\subsection{Definitions}

\begin{defn}
\label{partialisometries}
(a)  Let $K$ be a compact $\R$-tree. A {\em partial isometry} of $K$ is
an isometry between two closed subtrees of $K$. It is said to be
{\em non-empty} if its domain is non-empty.

\smallskip
\noindent
(b) A {\em system of 
isometries} $\CK=(K,\CA)$ consists of a compact
$\R$-tree $K$ and a finite set $\CA$ of non-empty partial isometries
of $K$.  This defines a {\em pseudo-group of partial isometries} of
$K$ by admitting inverses and composition. 
\end{defn}

We note that in the literature mentioned above it is usually required
that $K$ is a {finite tree}, i.e.  $K$ is a metric realisation of a
finite simplicial tree, or, equivalently, $K$ is the convex hull of
finitely many points. The novelty here is that we only require $K$ to
be compact.  Recall that a compact $\R$-tree $K$ may well have
infinitely many branch points, possibly with infinite valence, and
that $K$ may well contain finite trees of unbounded volume (but of
course $K$ has finite diameter).
In the context of this paper, however, all trees have a countable number of branch points, which makes compact trees slightly more tractable.

Any element of the free group $\FN$ over the basis $\CA$, given as
reduced word $w=z_1\ldots z_n \in \FA$, defines a (possibly empty)
partial isometry, also denoted by $w$, which is defined as the
composition of partial isometries $z_1 \circ z_{2} \circ \ldots \circ
z_n$. We write this {\it pseudo-action} of $\FA$ on $K$ on the right,
i.e.
$$x (u \circ v) = (xu)v$$
for all $x \in K$ and $u, v \in \FA$.  For any points $x,y\in K$ and
any $w \in \FA$ we obtain
$$x w=y$$
if and only if $x$ is in the domain $\dom(w)$ of $w$ and is sent by
$w$ to $y$.

A reduced word $w\in\FA$ is called {\em admissible} if it is non-empty as
a partial isometry of $K$.

\subsection{The $\R$-tree associated to a system of isometries}

A system of isometries $\CK=(K,\CA)$ defines an $\R$-tree $T_\CK$, 
provided with
an action of the free group 
$\FN = \FA$ by isometries. The construction
is the same as in the case where $K$ is a finite tree and will be 
recalled now.

As in \cite{gl} the tree $T_\CK$ can be described using a foliated
band-complex, but for non-finite $K$ one would not get a CW-complex.
We use the following equivalent construction in combinatorial terms.

The tree $T_\CK$ is obtained by gluing countably many copies of $K$
along the partial isometries, one for each element of $\FN$.  On the
topological space $\FN\times K$ these identifications are made formal
by defining
\[
T_\CK=\FN\times K/\sim 
\]
where the equivalence relation $\sim$ is defined by:
\[
(u, x)\sim (v, y) \iff x (u\inv v)=y
\]

The free group $\FN$ acts on $T_{\CK}$, from the left: this action is
simply given by left-multiplication on the first coordinate of each
pair $(u, x) \in \FN \times K$:
$$
w(u, x) = (wu, x)
$$
for all $u,w \in \FN, x \in K$.

Since $\FN$ is free over $\CA$, each copy $\{u\}\times K$ of $K$
embeds canonically into $T_\CK$.  Thus we can identify $K$ with the
image of $\{1\}\times K$ in $T_{\CK}$, so that every $\{u\} \times K$
maps bijectively onto $u K$.  Using these bijections, the metric on
$K$ defines canonically a pseudo-metric on $T_{\CK}$.  Again, by the
freeness of $\FN$ over $\CA$, this pseudo-metric is a metric.  The
arguments given in the proof of Theorem I.1. of \cite{gl} extend
directly from the case of finite $K$ to compact $K$, to show:

\begin{thm}
\label{thm:TK}
Given a system of isometries $\CK = (K, \CA)$ on a compact $\R$-tree
$K$, there exists a unique $\R$-tree $T_\CK$, provided with a left-action of $\FA$ by isometries, which satisfies:
\begin{enumerate}
\item[(1)] $T_\CK$ contains $K$ (as an isometrically embedded subtree).
\item[(2)] If $x \in K$ is in the domain of $a\in\CA$, then $a\inv x=xa$.
\item[(3)] Every orbit of the $\FA$-action on $T_{\CK}$ meets $K$.
  Indeed, every segment of $T_\CK$ is contained in a finite union of
  translates $w_{i} K$, for suitable $w_{i} \in\FA$.
\item[(4)] If $T$ is another $\R$-tree with an action of $\FA$ by
  isometries satisfying (1) and (2), then there exists a unique
  $\FA$-equivariant morphism $j: T_\CK\to T$ such that $j(x)=x$ for
  all $x\in K$.  \qed
\end{enumerate}
\end{thm}

\subsection{Systems of isometries
induced by an 
$\FN$-action on an $\R$-tree}
\label{inducedpartialisometries}

Frequent and 
important examples of systems of
isometries occur in the following context:

Let $T$ be any $\R$-tree with an $\FA$-action by isometries. Then 
any
 compact subtree $K \subset T$, which is sufficiently large
so that it intersects for any $a_{i} \in \CA$ the translate $a_{i} K$,
defines canonically a system of isometries given by:
$$
a_{i}:\begin{array}[t]{rcl} 
a_{i} K\cap K&\to& K \cap a_{i}\inv K\\
x&\mapsto &x a_{i} = a_{i}^{-1} x
\end{array}
$$
Since $K$ embeds into $T$, Theorem \ref{thm:TK} gives a map
$$
j: T_{\CK} \to T.
$$

The map $j$ fails in general to be injective.
 A classical technique for the study of an action on an $\R$-tree $T$
is to view $T_{\CK}$ as an approximation of $T$, and to consider a
sequence of increasing $K$. As $K$ increases to exhaust $T$, the
convergence of the sequence of $T_{\CK}$ to $T$ is well understood.
Moreover, if $K$ is a finite subtree of $T$, then $T_{\CK}$ is called
geometric and the full strength of the Rips machine can be used to
study it

In this article, we propose a new approach to study $T$, namely we
prove that there exists a compact subtree $K$ of $\bar T$ such that
$j$ is an isometry.  This gives the possibility to extend the results
proved for geometric trees (i.e. when $K$ is finite) to the case where
$K$ is only assumed to be compact.

\subsection{Basic lemmas}

We now present some basic lemmas about the action on $T_\CK$, for
admissible and non-admissible words in the given system of isometries.  
We first observe:

\begin{rem}
\label{disjointsubtrees}
(a)  Let $K$ and $K'$ be two
closed
disjoint subtrees of $T$. Then there exists a 
unique segment $[x, x']$ which {\em joins $K$ to $K'$}, i.e. one has 
$K \cap [x, x'] = \{x\}$ and $K' \cap [x, x'] = \{x'\}$. For any 
further points $y \in K, y' \in K'$  the segment $[y, y']$ contains 
both segments $[x, y']$ and $[x', y]$, and both contain $[x, x']$.

\smallskip
\noindent
(b)  As a shorthand, we use in the situation given above the following 
notation:
$$[K, K'] := [x, x'], \, \, \, [y, K'] := [y, x'], \, \, \, [K, y'] := [x, y'] $$

\smallskip
\noindent
(c)  If $y \in K$, then we set $[y, K] = [K, y] = \{y\}$, i.e. the 
segment of length 0 with $y$ as initial and terminal point.
\end{rem}

The following is a specification of statement (3) of Theorem \ref{thm:TK}:

\begin{lem}
\label{unioncontains}
For any non-admissible word $w\in\FA$ one has
\[
[K, w K] \,\,\,\subset \,\,\, \bigcup_{i=0}^{|w|}w_{i} K \, ,
\]
where $w_i$ is the prefix of $w$ with length $|w_{i}| = i$.
\end{lem}

\begin{proof}
  It suffices to show that for the reduced word $w = z_{1} \ldots
  z_{n}$ the union $\overset{n}{\underset{i\,=\,0}{\bigcup}} w_{i} K$
  is connected.  This follows directly from the fact that for all $i =
  1, \ldots, n$ the union $w_{i-1}K \cup w_{i}K = w_{i-1} (K \cup
  w_{i-1}\inv w_{i} K)$ is connected, since $w_{i-1}\inv w_{i} = z_{i}
  \in \CA^{\pm 1}$, and all partial isometries from $\CA$ are assumed
  to be non-empty.
\end{proof}

\begin{lem}\label{lem:admdom}
Let $\CK = (K, \CA), T_{\CK}$ and $\FA$ be as above.
\begin{enumerate}
\item[(1)] For all $w \in \FA$ one has
$$
\dom (w) =K\cap wK \, .
$$ 
\item[(2)] A word $w\in\FA$ is admissible if and only if $K\cap
  wK\neq\emptyset$.
\item[(3)] If $x \in \dom(w)$, then
$$w\inv x = x w \, .$$
\end{enumerate}
\end{lem}

\begin{proof}
Let $w \in \FA$ and 
$x \in T_{\CK}$.
If  
$x \in \dom(w) \subset K$, then 
the definition of $T_\CK$ gives $(1, x)\sim(w, xw)$, or equivalently
(compare Theorem \ref{thm:TK})
$$w\inv x = x w \, .$$
Therefore 
$x$
is contained in both $K$
and 
$wK$.
This shows:
$$\dom(w) \subset K \cap wK$$

Conversely, let $x$ be in $K\cap wK$. Then
$(1, x)\sim(w, y)$ for some point $y \in K$,
and by definition of $\sim$ the point $x$ lies
in the domain 
of $w$, with $xw=y$. Thus $w$ is admissible, and 
\[
K \cap wK \subset\dom(w) \, .
\]
\end{proof}

\begin{lem}
\label{prefixes}
For all $w \in \FA$
the following holds, where $w_{k}$ denotes the 
prefix of $w$ of length $k$:
\begin{enumerate}
\item[(1)]
$\quad \dom (w) \subset \dom (w_{k}) \quad$ for all $k \leq |w|$.   
\item[(2)]
$\quad \dom (w) =\overset{|w|}{\underset{k\,=\,0}{\bigcap}} w_{k} K
%\bigcap_{i=0}^{|w|}w_{i}K
$
\end{enumerate}
\end{lem}

\begin{proof}
  Assertion (1) follows directly from the definition of $\dom(w)$.
  Assertion (2) follows from assertion (1) and Lemma \ref{lem:admdom}
  (1).
\end{proof}

\begin{rem}
\label{leftversusright}
We would like to emphasize that it is important to keep the
$\FA$-action on $T_{\CK}$ apart from the $\FA$-pseudo-action on $K$.
This is the reason why we define the action on $T_{\CK}$ from the
left, whereas we define the pseudo-action by partial isometries on $K$
from the right.

This setting is also convenient to keep track of the two actions: a
point $x \in K$ lies in the domain of the partial isometry associated
to $w \in \FA$ if and only if $x$ is contained in $w K$ (Lemma
\ref{lem:admdom} (1)).  More to the point, the sequence of partial
isometries given by the word $w = z_{1} \ldots z_{n}$ defines points
$x z_{1} \ldots z_{i}$ which lie all inside of $K$ if and only if the sequence of
isometries of $T$ given by the prefixes of $w$ moves $K$ within $T$ in
such a way that $x$ is contained in each of the translates $z_{1}
\ldots z_{i}K$ (see Lemma \ref{prefixes} (2)).
\end{rem}

\begin{lem}\label{lem:lien}
  (a) For any non-admissible word $w\in\FA$ and any disjoint closed
  subtrees $K$ and $wK$, the arc $[K, wK]$ intersects all $w_iK$, where
  $w_i$ is a prefix of $w$.

\smallskip
\noindent
(b) For any point $Q \in K$ and any (possibly admissible) word $w \in
\FA$, the arc $[Q, wK]$ intersects all $w_{i} K$.
\end{lem}

\begin{proof}
(a) We prove part (a) by induction on the length of $w$. 

Let $u$ be the longest admissible prefix of $w$.  Thus $u\neq 1$, as
all partial isometries in $\CA^{\pm 1}$ are non-empty.  Hence we can
assume by induction that $u\inv w$ is either admissible or satisfies
the property stated in part (a).

Let $a$ be the next letter of $w$ after the prefix $u$. We write $w$
as reduced product $w=u \cdot a \cdot v$.  According to
Lemma~\ref{lem:admdom}~(2) one has:
\begin{enumerate}
\item[(i)]
$u K\cap K =\dom (u) \neq\emptyset$
\item[(ii)]
$u K \cap u a K =u \, \dom (a)\neq\emptyset$, \mbox{and }
\item[(iii)]
$ K \cap u a \, K =\emptyset$
\end{enumerate}
By (iii) there is a non-trivial segment $\beta = [K, uaK] \subset
T_{\CK}$ that intersects $K$ and $uaK$ only in its endpoints.  By (i)
and (ii) the segment $\beta$ is contained in the subtree $uK$: there
are points $x, y \in K$ such that $\beta = [u x, uy]$.  Since $u x$
belongs to $K \cap uK = \dom(u)$, it follows from Lemma \ref{prefixes}
(2) that $ux$ also belongs to every $u' K$, for any prefix $u'$ of
$u$.

Moreover, for any prefix $v'$ of $v$ one has, by 
Lemma~\ref{lem:admdom}~(1) and Lemma~\ref{prefixes}~(1):
\[
uav' K\cap u K \, \, =\, \,  u \, \dom (a v') \, \, \, \subset \, \, \,
u \, \dom
(a) \, \, = \, \,  ua K \cap u K
\]
From this we deduce that
\[
\begin{array}[t]{llllll}
uav' K \cap [ux,uy] & \subset & uav' K \cap [ux,uy]  \cap u K & &
 \\
& \subset &
[ux, uy] \cap ua K \cap u K &&&\\
&\subset &[ux, uy] \cap ua K & =&\{uy\} \, .
\end{array}
\]
Since the segment $\alpha = [K, w K]$ is by Lemma \ref{unioncontains}
contained in the union
\[
\bigcup_{i=0}^{|w|}w_{i} K
\]
it follows from the above derived inclusion $uav' K \cap [ux,uy]
\subset \{uy\}$ that $\alpha$ is the union of $\beta = [ux, uy]$ and
of the segment $\gamma = [uy, w K]$, with $\beta \cap \gamma =
\{uy\}$.

If $av$ is admissible, then the endpoint of $\gamma$ is contained in
the intersection of all $uav'K$, by Lemma \ref{prefixes} (2).  If $av$
is non-admisible, we apply the induction hypothesis to $u\inv w = av$
and obtain that every $av' K$ meets the arc $\gamma' = [K, avK]$. But
$u \gamma'$ is a subarc of $\gamma$, so that the arc $[ux, uy] \cup
\gamma$ meets infact all $w_{i} K$, as claimed.

\smallskip
\noindent
(b) In case that $w$ is non-admissible, there is a largest index $i$
such that $K \cap w_{i}K \neq \emptyset$.  We can now apply statement
(a) to $w_{i}\inv K$ and $w_{i}\inv w$ to get the desired conclusion.

If $w$ is admissible, then $\dom(w) = K \cap wK$ (by Lemma
\ref{lem:admdom} (1)). Hence the arc $[Q, wK]$ is contained in $K$,
and by Lemma \ref{prefixes} (2) its endpoint is contained in any
$w_{i} K$.
\end{proof}

\begin{lem}
\label{center}
Let $w, w' \in \FA$ with maximal common prefix $u \in \FA$.
Then for any triplet of points $Q \in K, \, R \in wK$ and $R' \in w' K$
the arcs $[Q, R]$ 
and $[Q, R']$ 
intersect in an arc $[Q, P]$ with endpoint $P \in u K$.
\end{lem}

\begin{proof}
Let $[Q, Q_{1}]$ the arc which joins $K$ to $uK$.  It follows directly 
from 
Lemma \ref{lem:lien} (b)
that $Q_{1}$ lies on both $[Q, R]$ and $[Q, 
R']$. Similarly, let $[R, R_{1}]$ and $[R', R'_{1}]$ be the arcs that 
join $R$ to $uK$ and $R'$ to $uK$ respectively. After applying $w\inv$ 
or ${w'}\inv$ we obtain in the same way that 
$R_{1}$ lies on both $[Q, R]$ and $[R, R']$, and that
$R'_{1}$ lies on both $[Q, R']$ and $[R, R']$. Hence the geodesic 
triangle in $T_{\CK}$ with endpoints $Q, R, R'$ contains the geodesic 
triangle with endpoints $Q_{1}, R_{1}$ and $R'_{1}$, and the center of 
the latter is equal to the center $P$ of the former.  But $Q_{1}, R_{1}$ and 
$R'_{1}$ are all three contained in $uK$, so that $P$ is contained in 
$uK$.
\end{proof}

In the following statement and its proof we use the standard 
terminology for group elements acting on trees, as recalled in \S\ref{prelims} above.

\begin{prop}\label{prop:axes}
Let $w \in \FA$ is any cyclically reduced word. If the action of $w$
on $T_{\CK}$ is hyperbolic,
then the axis of $w$ intersects $K$.
If the action of $w$ on $T_{\CK}$
is elliptic, then $w$ has a fixed point in $K$.
\end{prop}

\begin{proof}
  If $w$ is not admissible, let $[x, wy]$ be the segment that joins
  $K$ to $w K$: these two translates are disjoint by Lemma
  \ref{lem:admdom} (2). As $w$ acts as an isometry, $[wx, w^{2} y]$ is
  the segment that joins $w K$ to $w^{2} K$. Moreover, since $w$ is
  assumed to be cyclically reduced, the segment that joins $K$ to
  $w^{2} K$ intersects $w K$, by Lemma \ref{lem:lien}.

  Any two consecutive segments among $[x, wy]$, $[wy, wx]$, $[wx,
  w^{2}y]$ and $[w^{2}y, w^{2} x]$ have precisely one point in common,
  by Remark \ref{disjointsubtrees}, and hence their union is a
  segment. This proves that $ w x$ belongs to $[x, w^{2} x]$, and that
  $x$ is contained in the axis of $w$.

  If $w$ is admissible, then either there exists $n \geq 0$ such that
  $w^n$ is not admissible, in which case we can fall back on the above
  treated case, as $w$ and $w^n$ have the same axis.  Otherwise, for
  arbitrary large $n$ there exists a point $x \in K$ such that $w^n x
  \in K$, by Lemma \ref{lem:admdom} (2).  But $K$ is compact and hence
  has finite diameter.  This implies that the action of $w$ on $T$ is
  not hyperbolic, and hence it is elliptic: $w$ fixes a point of $T$.
  Some such fixed point lies on $[x, wx]$ (namely its center), and
  hence in the subtree $K$.
\end{proof}

\subsection{Admissible laminations}
\label{admissiblelamination}

In this subsection we use the concepts of {\em algebraic lamination},
{\em symbolic lamination} and {\em laminary language} as defined in
\cite{chl1-I}, and the equivalence between these three points of view
shown there.  The definitions and the notation have been reviewed in \S\ref{algebraiclaminations} above.

For any system of isometries $\CK=(K,\CA)$ denote by $\Adm(\CK)\subset
\FA$ the set of admissible words.  The set $\Adm(\CK)$ is stable with
respect to passage to subwords, but it is not {\em laminary} (see
\cite[Definition 5.2]{chl1-I}): not every admissible word $w$ is
necessarily equal, for all $k \in \N$, to the word $v\chop_{k}$
obtained from some larger $v \in \Adm(\CK)$ by ``chopping off'' the
two boundary subwords of length $k$.  As does any infinite subset of
$\FA$, the set $\Adm(\CK)$ generates a laminary language, denoted
$\CLadm(\CK)$, which is the largest laminary language made of
admissible words:
$$\CLadm(\CK) = \{ w \in \FA \mid \forall k\in\N \,\, \exists v\in
\Adm(\CK): 
w = v\chop _{k}\}$$
Clearly one has $\CLadm(\CK) \subset \Adm(\CK)$, but the converse is 
in general false.

As explained in
\S\ref{algebraiclaminations},
any laminary language determines an algebraic lamination (i.e. a 
closed $\FN$-invariant and flip-invarinat subset of $\partial^{2}\FN$), and 
conversely. The algebraic lamination determined by $\CLadm(\CK)$ is 
called {\em admissible lamination}, and denoted by $\Ladm(\CK)$.

\smallskip

An infinite word $X\in\partial\FA$ is {\em admissible} if all of its
prefixes $X_{n}$ are admissible.  The set of admissible infinite words
is denoted by $\Loneadm (\CK)$. It is a closed subset of $\partial\FA$
but it is not invariant under the action of $\FA$.

For any infinite admissible $X$ the {\em domain} $\dom(X)$ of $X$ is
defined to be the intersection of all domains $\dom(X_{n})$.  Since
$K$ is compact, one has
$$
\dom(X) \neq \emptyset
$$
for all $X \in \Loneadm (\CK)$.

\smallskip

A biinfinite indexed reduced word $Z = \ldots z_{-1} z_{0} z_{1}
\ldots$, with $z_{i} \in \CA^\pm$, is called {\em admissible}, if its
two halves $Z^+ = z_{1} z_{2} \ldots$ and $Z^- = z_{0}\inv z_{-1}\inv
\ldots$ are admissible, and if the intersection of the domains of
$Z^+$ and $Z^-$ is non-empty. The {\em domain} of $Z$ is defined to be
this intersection:
$$
\dom(Z) = \dom(Z^{+}) \cap \dom(Z^-)
$$
We observe that $Z$ is admissible if and only if all its subwords are
admissible.

The set  
of biinfinite admissible words is called the {\em admissible symbolic
lamination} of the system of isometries $\CK = (K,\CA)$.

\medskip

We use now the notion of the dual lamination of an $\R$-tree with 
isometric $\FN$-action as introduced in \cite{chl1-II} and reviewed 
above in \S\ref{duallamination}.

\begin{prop}\label{prop:dualvsadm}
For any system of isometries $\CK$
one has
\[
L(T_\CK)\subseteq\Ladm(\CK).
\]
\end{prop}

\begin{proof}
  Let $u \in \FA$ be a non-admissible word, and let $\varepsilon=d(K,
  u K)$. By Lemma \ref{lem:admdom} (2) one has $\varepsilon>0$. Let
  $w$ be a cyclically reduced word that contains $u$ as a subword: we
  write $w=u_1 \cdot u \cdot u_2$ as a reduced product.  By
  Proposition~\ref{prop:axes}, the axis of $w$ passes through $K$. But
  if $x$ is any point in $K$, the segment $[ x, w x]$ contains the
  segment that joins the disjoint subtrees $u_{1} K$ and $u_{1} u K$,
  by Lemma \ref{lem:lien}, and hence the translation length of $w$,
  which is realized on its axis, is bigger than $\varepsilon$. This
  proves that $u$ is not in $\L(T_\CK)$ (see \S\ref{duallamination}) .

  As the laminary language of $\CLadm(\CK)$ is the largest laminary
  language made of admissible words, this concludes the proof.
\end{proof}

\section{The map $\CQK$ for a system of isometries}
\label{sectionfour}

In this section we define the map $\CQK$ and we prove that it  is the equivalent of the
map $\CQ$ from \S\ref{sec:Q}, for systems of isometries $\CK$. For this definition we distinguish
two cases: If $X\in\partial\FA$ is not eventually admissible we define
$\CQK(X)$ in \S\ref{sec:rays}. If $X$ is eventually admissible, the
definition of $\CQK(X)$ is given in \S\ref{sec:QKadm}, and in this case we need the hypothesis  that
the system of isometries has independent generators.  Both cases are
collected together in \S\ref{sec:QK} to obtain a continuous equivariant map $\CQK$.

\subsection{The map $\CQK$ for non-eventually admissible words}\label{sec:rays}

As in \S\ref{sectionthree}, let $\CK=(K,\CA)$ be a system of
isometries on a compact $\R$-tree $K$, and let $T_\CK$ be the
associated $\R$-tree, provided with an action of the free group $\FA$
by isometries. Let $X\in\partial\FA$ be an infinite reduced word and
denote as before by $X_i$ the prefix of $X$ of length $i\geq 0$.

\begin{defn}\label{def:nonU}
  An infinite word $X\in\partial\FA$ is {\em eventually admissible} if
  there exists an index $i$ such that the suffix ${X_i}\inv X$ of $X$
  is admissible.
\end{defn}

Note that an infinite word $X\in\partial\FA$ is not eventually admissible if for every index $i
\geq 0$ there is an index $j>i$ such that the subword $X_{[i+1,j]}=X_i\inv X_j$
of $X$ between the indices $i+1$ and $j$ is not admissible.

Let $X\in\partial\FA$ be not eventually admissible, and let
$i_0>0$ be such that the prefix $X_{i_0}$ of $X$ of length $i_0$ is
not admissible.  Then for any $i\geq i_0$, the prefix $X_i$ is not admissible, and
thus, by Lemma~\ref{lem:admdom}, $K$ and $X_iK$ are disjoint. By
Lemmas~\ref{lem:lien} and \ref{center}, for any $j\geq i\geq i_0$ the segment
$[K,X_iK]$ and $[K,X_jK]$ are nested and have the same initial point
$Q\in K$.  Let $Q_i$ be the terminal point of $[K,X_iK]$:
$$
[Q,Q_i]=[K,X_iK]
$$ The sequence of $Q_i$ converges in $\hat T_\CK$ with respect to
both the metric and the observers' topology. Moreover, the two limits
are the same.

\begin{defn}\label{def:QKnonU}
For any  $X\in\partial\FA$ which is  not eventually admissible, we define:
$$
\CQK(X)=\lim_{i\to\infty} Q_i
$$ 
\end{defn}

\begin{prop}
\label{prop:QKnonU}
Let $\CK=(K,\CA)$ be a system of isometries on a compact $\R$-tree $K$.
Let $X \in \partial \FA$ be  not eventually admissible.

Let $w_{n} \in \FA$ be a sequence of words which converge in $\FA \cup
\partial \FA$ to $X$, and let $P_n\in w_nK$.  Then the sequence of
points $P_{n}$ converges in $\Tobs_{\CK}$ to $\CQ_\CK(X)$, and
$\CQ_\CK(X)$ belongs to $\Tobs_{\CK} \smallsetminus T_{\CK}$.
\end{prop}
\begin{proof}
  We use the above notations.  For every index
  $i\geq 0$, let $[Q_i,R_i]$ be the intersection of $[Q,\CQ_\CK(X)]$
  with $X_iK$. Hence for $i\geq i_0$ the point $Q_i$ is, as before, the terminal
  point of the segment $[K,X_iK]$. The segments $[Q,Q_i]$ are increasingly nested,
  the segments $[R_i,\CQK(X)]$ are decreasingly nested, $Q_i$ is a
  point of $[Q,R_i]$ and $R_i$ is a point of $[Q_i,\CQK(X)]$.

As $X$ is not eventually admissible, for every
  index $i \geq 0$ there is an index $j > i$ such that the subword
  $X_{[i+1,j]}$ of $X$ between the indices $i+1$ and $j$ is not
  admissible. By Lemma~\ref{lem:admdom} the segments $[Q_{i}, R_{i}]$
  and $[Q_{j}, R_{j}]$ are disjoint.

  For any $n$, let $i(n)$ be the length of the maximal common prefix
  of $w_n$ and $X$.  By Lemma~\ref{center}, the maximal common segment
  $[Q,P'_n]$ of $[Q,P_n]$ and $[Q,\CQ_\CK(X)]$ has its terminal point
  $P'_n$ in $[Q_{i(n)},R_{i(n)}]$.  As $X$ is not eventually
  admissible, for $m$ big enough the subword $X_{[i(n)+1,i(m)]}$ of $X$
  between the indices $i(n)+1$ and $i(m)$ is not admissible and the
  segments $[Q_{i(n)},R_{i(n)}]$ and $[Q_{i(m)},R_{i(m)}]$ are
  disjoint. Therefore the maximal common segment of $[Q,P_n]$ and
  $[Q,P_m]$ is also the maximal common segment of $[Q,P_n]$ and
  $[Q,\CQ_\CK(X)]$, and hence it is equal  to $[Q,P'_n]$.

  The points  $P'_n$ converge to $\CQK(X)$, as any sequence of points in
  $[Q_{i(n)},R_{i(n)}]$ does, and this proves that
\[
\liminf{}_Q P_n=\CQ_\CK(X).
\] 
By Lemma~\ref{lem:convergence} any subsequence of $P_n$, which
converges in $\Tobs_\CK$, necessarily converges to $\CQ_\CK(X)$. Hence
by compactness of $\Tobs_\CK$, the sequence of all of the points $P_n$
converges  to $\CQ_\CK(X)$ with respect to the observers' topology.

If $P$ is a point in $uK$ for some $u$ in $\FN$, then
the maximal common segment $[Q,P']$ of $[Q,P]$ and $[Q,\CQ_\CK(X)]$
has its endpoint $P'$ in $[Q_i,R_i]$, where $X_i$ is the maximal common
prefix of $u$ and $X$. Thus $P'\neq \CQK(X)$, and  hence $\CQ_\CK(X)$ is not contained in $T_\CK$.
\end{proof}

\subsection{Independent generators}
\label{independent-generators}

The following concept is due to Gaboriau \cite{gab}, in the case of
finite $K$, and we extend it here to the compact case.

\begin{defn}
\label{independentgenerators}
Let $\CK=(K,\CA)$ be a system of isometries on a compact $\R$-tree
$K$. Then $\CK$ is said to have {\em independent generators} if, for
any infinite admissible word $X\in\partial\FA$, the non-empty domain
of $X$ consists of exactly one point.
\end{defn}

The same arguments as in \cite{gab} show the following equivalences.
However, they will not be used in the sequel.

\begin{rem}\label{rem:indgen}
Let $\CK=(K,\CA)$ be a system of isometries on a compact $\R$-tree
$K$. The following are equivalent:
\begin{enumerate}
\item[(1)] $\CK$ 
has independent generators.
\item[(2)] 
Every non-trivial admissible word
fixes at most one point of $K$.
\item[(3)] The action of $\FA$ on the associated tree $T_\CK$ has trivial
arc stabilizers.
\end{enumerate}
\end{rem}

Note that Gaboriau \cite{gab} used originally property (2) as 
definition, but in our context this seems less natural.

\subsection{The map $\CQK$ for eventually admissible words}\label{sec:QKadm}

Let $\CK=(K,\CA)$ be a system of isometries on a compact $\R$-tree
$K$.  Consider the set $\Loneadm(\CK)\subset\partial\FA$ of infinite
admissible words as defined in \S\ref{admissiblelamination}.  

\begin{defn}
\label{def:QKadm}
Let $\CK$ be a system of isometries which
has independent generators.
Then for any
infinite admissible word $X\in\Loneadm(\CK)$ there exists exactly one
element of $K$ in the domain of $X$, which will be called 
$\CQ_{\CK}(X)$. 
\end{defn}

\begin{lem}
\label{prefixequivariance}
Identify $K$ with the image of $\{1\} \times K$ in $\hat T_{\CK}$ as
in \S\ref{sectionthree}, and let $X \in \Loneadm(\CK)$.
\begin{enumerate}
\item[(1)]
Denoting as before by $X_{i}$ the prefix of $X$ of length $i \geq 1$, 
we obtain:
$$\{\CQK(X)\} = 
\bigcap_{i \, \geq \, 1} \, X_{i} K$$
\item[(2)]
For every $i \geq 1$ we have:
$$\CQK(X_{i}\inv X) = X_{i}\inv \CQK(X)$$
\end{enumerate}
\end{lem}

\begin{proof}
Assertion (1) follows directly from Lemma \ref{prefixes} (2) and the 
above definition of the map $\CQK$.
Assertion (2) follows directly from (1).
\end{proof}

Recall from Definition~\ref{def:nonU} that an infinite words
$X\in\partial\FA$ is eventually admissible if it has a prefix $X_{i}$
such that the infinite remainder $X'_{i} = X_{i}\inv X$ is admissible.
We observe that for all integers $j \geq i$ the word $X_{i}\inv X_{j}$
is admissible, so that Lemma \ref{prefixequivariance} (2) gives:
$$X_{i} \CQ_{\CK}(X'_{i}) 
= X_{i} \CQ_{\CK}(X_{i}\inv X_{j} X'_{j})
= X_{i} (X_{i}\inv X_{j}) \CQ_{\CK} (X'_{j})
=  X_{j} \CQ_{\CK} (X'_{j})
$$
Hence the following definition does not depend on the choice of the 
index $i$.

\begin{defn}\label{def:QKU}
For any eventually admissible word $X\in\partial\FA$ we define
$$
\CQ_{\CK}(X) = X_{i} \CQ_{\CK}(X'_{i}).
$$
\end{defn}

We note that for any element $u\in\FA$ and any eventually admissible word  $X\in\partial\FA$ one has:
$$
\CQK(uX)=u\CQK(X)
$$

\begin{prop}\label{prop:QKU}
  Let $\CK=(K,\CA)$ be a system of isometries on a compact $\R$-tree
  $K$ with indenpendent generators.  Let $X\in\partial\FA$ be an
  eventually admissible word.

  For any element $P$ in $T_\CK$, and any sequence $w_n$ of elements
  of $\FN$ that converge to $X$, the sequence of points $w_nP$
  converges to $\CQ_\CK(X)$, with respect to the observers' topology
  on $T_\CK$.
\end{prop}
\begin{proof}
  Up to multiplying by the inverse of a prefix we can assume that $X$
  is admissible and $\CQK(X)\in K$. By compactness of $\Tobs_\CK$ we
  can assume that $w_nP$ converges to some point $Q_\infty$.  By contradiction
  assume that $Q_\infty\neq\CQ_\CK(X)$, and let $M$ be a point in the
  open interval $(Q_\infty,\CQK(X))$. From Lemma~\ref{lem:convergence} we deduce 
\[
Q_\infty=\liminf{}_{\CQK(X)} w_nP.
\]
Thus, for $n$ and $m$ big enough, the maximal common segment
$[\CQK(X),P_{m,n}]$ of the segments $[\CQK(X),w_nP]$ and $[\CQK(X),w_mP]$ contains
$M$. As $w_n$ converges to $X$, for $n$ fixed and for $m$ sufficiently
large, the maximal common prefix of $w_n$ and $w_m$ is a
prefix $X_i$ of $X$. By Lemma~\ref{center}, $P_{m,n}$ is contained in $X_iK$.
By Lemma~\ref{prefixequivariance}, $\CQK(X)$ is also contained $X_iK$,
and hence, so is $M$. As $m$ and $n$ grow larger, the index $i$ goes to
infinity (since $w_n\to X$), which proves that $M$ is contained in the
intersection of all the $X_iK$. Since we assumed $M\neq\CQK(X)$, this
contradicts the independent generators' hypothesis.
\end{proof}

\subsection{Continuity of the map $\CQ_\CK$}\label{sec:QK}

As any element of $\partial\FA$ is either eventually admissible or
not, from Definitions~\ref{def:QKnonU} and \ref{def:QKU} we collect a
map $\CQK$.

\begin{cor}\label{cor:QK}
  Let $\CK=(K,\CA)$ be a system of isometries on a compact $\R$-tree
  $K$ with independent generators.  The map
  $\CQ_\CK:\partial\FN\to\Tobs_\CK$ is equivariant and continuous.

  For any point $P$ in $T_\CK$, the map $\CQK$ defines the continuous
  extension to $\FN\cup\partial\FN$ of the map
\[
Q_P:\begin{array}[t]{rcl}
\FN&\to&\Tobs_\CK\\
w&\mapsto&wP
\end{array}
\]
\end{cor}
\begin{proof}
  Equivariance and continuity of $\CQK$ follow from the second part of
  the statement, which is proved in Propositions~\ref{prop:QKnonU}
  and \ref{prop:QKU}.
%, blended with a small dose of Bourbaki
%  extract\footnote{at this point the second cook ran out of the
%    kitchen, shouting undistinguishible words in french (``... farouche
%    ...'');} .
\end{proof}

\section{Proof of the Main Theorem}
\label{sectionfive}

Throughout this section let $T$ be an $\R$-tree provided with a
minimal, very small action of $\FN$ by isometries which has dense
orbits.  Hence we obtain from Theorem~\ref{thm:Qexists} 
an equivariant and continuous map $\CQ$, which we denote here by
$\CQT:\partial\FN\to\Tobs$.

Let $\CA$ be a basis of $\FN$, and
let $K$ be a compact subtree of $\bar T$. Let $\CK=(K,\CA)$ be the induced
system of isometries $a_{i}: K \cap a_{i} K \to a_{i}\inv K \cap K, \,
x \mapsto x a_{i} = a_{i}\inv x$, as discussed in
\S\ref{inducedpartialisometries}.  We assume that $K$ is chosen large
enough so that for each $a_{i} \in \CA$ the intersection $K \cap a_{i}
K$ and hence the partial isometry $a_{i}\in\CA$ is non-empty.  As a
consequence (see \S\ref{sectionthree}), there exists an $\R$-tree
$T_\CK$ with isometric action by $\FN$, and by Theorem~\ref{thm:TK}
there exists a unique continuous $\FN$-equivariant map
$$
j: T_\CK\to \bar T 
$$
which induces the identity map $T_{K} \supset K \overset{j}{\to} K 
\subset\bar T$.

\begin{lem}\label{lem:Kindependant}
  The system of isometries $\CK=(K,\CA)$
  has independent generators.
\end{lem}
\begin{proof}
  Let $Q$ be a point in the domain of an infinite admissible word $X$,
  compare \S\ref{admissiblelamination}.  Then for any prefix $X_n$ of
  $X$, the point $QX_n={X_n}\inv Q$ is also contained in $K$ (recall
  that we write the action of $\FA$ on $T_{\CK}$ on the left, and the
  pseudo-action of partial isometries of $\CK$ on the right).   

  By Theorem~\ref{thm:TK}, $j$ restricts to an isometry between
  $K\subset T_\CK$ and $K\subset T$. Therefore, for any $n\geq 0$,
  ${X_n}\inv j(Q)$ lies in $K\subset T$.  By Lemma~\ref{lem:QTandK}, we get
  $\CQT(X)=j(Q)$.

  This proves that the domain of $X$ consists of at most the point
  $j\inv(\CQT(X))$.  Hence $\CK$ has independent generators.
\end{proof}

As a consequence of Lemma~\ref{lem:Kindependant}, we can apply Corollary~\ref{cor:QK} to obtain an equivariant and
continuous map $\CQK:\partial\FN\to\Tobs_\CK$.

\begin{lem}\label{lem:jQX}
For any $X\in\partial\FN$ such that $\CQK(X)$ is contained in $T_\CK$, one has
\[
j(\CQK(X))=\CQT(X).
\]
\end{lem}
\begin{proof}
  By Proposition~\ref{prop:QKnonU}, $X$ is eventually
  admissible and by equivariance of $\CQK$, $\CQT$ and $j$, we can
  assume that $X$ is admissible and that $\CQK(X)$ is in $K$. By
  Definition~\ref{def:QKadm}, for any $i\geq 0$,
  $\CQK(X)\cdot X_i={X_i}\inv\CQK(X)$ lies in $K$. 
  
  By Theorem~\ref{thm:TK}, $j$ restricts to an isometry between
  $K\subset T_\CK$ and $K\subset T$. Therefore for any $i\geq 0$, the
  point ${X_i}\inv j(\CQK(X))$ lies in $K\subset T$. Thus we can apply
  Lemma~\ref{lem:QTandK} to get $\CQT(X)=j(\CQK(X))$.
\end{proof}

\begin{lem}\label{lem:lamsubset}
The admissible lamination of $\CK$ is contained in the dual lamination of $T$:
$$\Ladm(\CK)\subset L(T)$$
\end{lem}
\begin{proof}
  The admissible lamination $\Ladm(\CK)$ (see
  \S\ref{admissiblelamination}) is defined by all biinfinite words
  $Z$ in $\CA^\pm$ such the two half-words $Z^{+}$ and $Z^{-}$ have
  non-empty domain, and the two domains intersect non-trivially.  Thus
  $\CQK(Z^{+}) = \CQK(Z^{-})$ is a point in $K$. Thus by
  Lemma~\ref{lem:jQX} one has $\CQT(Z^{+}) = \CQT(Z^{-})$. The latter
  implies (and is equivalent to) that $Z$ belongs to $L(T)$.
\end{proof}

We sumarize the above discussion in the following commutative diagram:
 \[
 \xymatrix{
 &\partial F_N\ar[ld]_\CQK \ar@{>>}[rd]^\CQT\\
 \Tobs_{\CK}&&\Tobs\\
%\rule{0cm}{.5cm}
T_\CK \ar@{^{(}->}[u]\ar[rr]^j&&
%\rule{0cm}{.5cm}
\bar T \ar@{^{(}->}[u]
}
\]
All the maps in the diagram are equivariant and continuous, where the
topology considered on the bottom line is the metric topology.

We can now prove the main result of this paper.  Recall from
\S\ref{limitset} that for any basis $\CA$ of $\FN$ and $T$ as above
the set $\Omega_{\CA} \subset \bar T$ denotes the limit set of $T$ with respect
to $\CA$.

\begin{thm}\label{thm:mainthm}
  Let $T$ be an $\R$-tree with very small minimal $\FN$-action by
  isometries, and with dense orbits. Let $\CA$ be a basis of $\FN$,
  and let $K \subset \bar T$ be a compact subtree which satisfies $K
  \cap a_{i} K \neq \emptyset$ for all $a_{i} \in \CA$.  Then the
  following are equivalent:

\begin{enumerate}
\item[(1)]
\label{jisometry} 
The restriction of the canonical map $j: T_{\CK} \to \bar T$ to the minimal 
$\FN$-invariant subtree $T_{\CK}^{min}$ of $T_{\CK}$
defines an isometry $j^{min}: T_{\CK}^{min} \to T$.
\item[(2)]
\label{laminationequal}
$L(T) \subset \Ladm(\CK)$  \qquad ($\, \Longleftrightarrow 
L(T)=\Ladm(\CK)\, ,$ by Lemma~\ref{lem:lamsubset})
\item[(3)]
\label{containsheart} 
$\Omega_\CA\subset K$
\end{enumerate}
\end{thm}

\begin{proof}
\underline{(1) $\Longrightarrow$ (2)}:
By the assumption on $j$
the minimal subtree $T_{\CK}^{min} \subset T_{\CK}$ is 
isometric 
to $T$. Hence the dual laminations satisfy $L(T) = 
L(T_{\CK}^{min})$, and by 
Remark \ref{dualgivenbyminimal}
one has $L(T_{\CK}^{min}) = L(T_{\CK})$. 
We now apply
Proposition
\ref{prop:dualvsadm} 
to get $L(T_{\CK})\subset\Ladm(\CK)$.

\medskip
\noindent
\underline{(2) $\Longrightarrow$ (3)}: By Definition
\ref{limitset-heart}, a point $Q \in T$ belongs to the limit set
$\Omega_{\CA}$ if and only if there is a pair of infinite words $(X,
Y) \in L^{2}(T) \subset \partial^{2} \FA$, with initial letters $X_{1}
\neq Y_{1}$, which satisfy $Q_{T}(X) = Q_{T}(Y) = Q$.  By assumption,
$L(T)$ is a subset of $\Ladm(\CK)$, so that the reduced words $X$, $Y$
and $X\inv \cdot Y$ are admissible for the system of isometries $\CK$.
By Definition~\ref{def:QKadm}, $\{\CQK(X)\}$ is the domain of $X$ and
$Y$, and thus is contained in $K$. We deduce from Lemma~\ref{lem:jQX} that 
$j(\CQK(X))=\CQT(X)=Q$, and $Q$ lies in $K$.

\medskip
\noindent
\underline{(3) $\Longrightarrow$ (2)}: Let $Z$ be a biinfinite indexed
reduced word in the symbolic lamination $L_\CA(T)$ defined by the dual lamination
$L(T)$ of $T$ (see \S\ref{algebraiclaminations}). That is to say,
$Z=(Z^-)\inv \cdot Z^+$, written as a reduced product, and
$\CQ_T(Z^-)=\CQ_T(Z^+)$ is a point $Q \in \Omega_\CA$. For any
$n\in\Z$, we consider the shift $\sigma^n(Z)$ of $Z$ as in
Remark~\ref{symbolicpartialisom}.  If $u$ is the prefix of $Z^+$ of
length $n$ (or, if $n<0$, the prefix of $Z^-$ of length $-n$), then
$\sigma^n(Z)= (Z^-)\inv u \cdot u\inv Z^+$ and $\CQ_T(u\inv
Z^+)=\CQ_T(u\inv Z^-) = u\inv Q$, and this is again a point of
$\Omega_\CA$ and thus contained in $K$, by hypothesis. Therefore, both
$Z^{+}$ and $Z^-$ are admissible, and $\dom(Z^{+}) = \dom(Z^-) =
\{Q\}$. Thus $Z$ is an admissible biinfinite word of the system of
isometries $\CK=(K,\CA)$, which shows $L(T) \subset \Ladm(\CK)$.

\medskip
\noindent
\underline{(2) $\Longrightarrow$ (1)}:
Since the dual lamination $L(T)$ is a subset of the
admissible lamination $\Ladm(\CK)$, for any pair of distinct infinite
words $X,Y\in\partial\FA$ the equality $\CQT(X)=\CQT(Y)$ implies that
$X\inv Y$ is admissible, and from Definition~\ref{def:QKadm} we deduce
$\CQK(X)=\CQK(Y)$.  Thus  the map $\CQK: \partial
\FN \to \Tobs_{\CK}$ factors over the quotient map $\pi:\partial \FN
\to \partial \FN / L^{2}(T)$ (see \S\ref{themapQ2}) to define an
equivariant map $s: \partial\FN/L^{2}(T)\to\Tobs_\CK$.

As the topology on $\partial\FN/L^{2}(T)$ is the quotient topology (see \S\ref{themapQ2}) and
as $\CQK$ is continuous (see Corollary~\ref{cor:QK}), the map $s$ is
continuous. Since $\varphi:\partial\FN/L^{2}(T)\to\Tobs$ is a
homeomorphism (see Theorem~\ref{thm:homeo}), we deduce that the image
of $s$ is an $\FN$-invariant connected subtree of
$\Tobs_\CK$. Therefore the image of $s$ contains the minimal subtree
$\TKmin$ of $T_\CK$.

As a consequence, for any point $P$ in $\TKmin$ there exists an
element $X\in\partial\FN$ such that $s(\pi(X))=\CQK(X)=P$. From
Lemma~\ref{lem:jQX} we obtain $\jmin(P)=j(\CQK(X))=\CQT(X)$. By
definition of the homeomorphism $\varphi$, one has
$\varphi\inv(\jmin(P))=\pi(X)$ and $s(\varphi\inv(\jmin(P)))=P$. This
proves that $\jmin$ is injective.

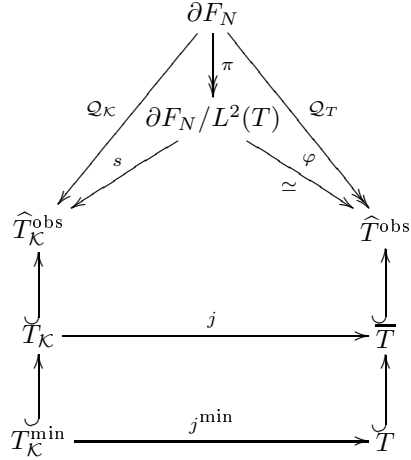
\[
 \xymatrix{
 &\partial F_N\ar[ldd]_\CQK \ar@{>>}[rdd]^\CQT \ar@{>>}[d]^\pi\\
&\partial F_N/L^2(T)\ar[ld]_s \ar[rd]^\varphi_\simeq\\
 \Tobs_{\CK}&&\Tobs\\
%\rule{0cm}{.5cm}
T_\CK \ar@{^{(}->}[u]\ar[rr]^j&&
%\rule{0cm}{.5cm}
\bar T \ar@{^{(}->}[u]\\
%\rule{0cm}{.5cm}
\TKmin \ar@{^{(}->}[u] \ar[rr]^\jmin&&
%\rule{0cm}{.5cm}
T \ar@{^{(}->}[u]
}
\]

Since $j$ is continuous with respect to the metric topology, since $j$
maps $K$ isometrically, and since $T_{\CK} = \FN K$, this implies that
$\jmin$ is an isometry.
${}^{}$\end{proof}

\bigskip

Recall from \S\ref{limitset} that the heart $K_{\CA} \subset \bar T$
denotes the convex hull of the limit set $\Omega_{\CA}$ of $T$ with
respect to the basis $\CA$.  We denote by $\CK_{\CA} = (K_{\CA}, \CA)$
the associated system of partial isometries.

We remark that, in the above theorem, the map $\CQK$ may fail to be
surjective onto $T_\CK$ if $K$ is too large. And hence, $j$ may fail
to be injective even if the limit set $\Omega_\CA$ is contained in
$K$. This is the reason why we considered the minimal subtree $\TKmin$ of
$T_\CK$. However if $K$ is exactly equal to the heart $K_\CA$ we get
the following corollary.

\begin{cor}\label{cor:KA}
  Let $T$ be an $\R$-tree with very small minimal $\FN$-action by
  isometries, and with dense orbits. Let $\CA$ be a basis of $\FN$,
  with heart $K_{\CA}$.  The map $j:T_{\CK_{\CA}}\to\bar T$ is
  isometric and its image contains $T$.
\end{cor}

\begin{proof}
  By definition, for $K = K_{\CA}$ the three equivalent conditions of
  Theorem~\ref{thm:mainthm} are satisfied.

  In the proof of implication (2)$\Rightarrow$(3) of
  Theorem~\ref{thm:mainthm}, we proved that $\Omega_\CA$ is in the
  image of $\CQK$.  In the proof of implication (2)$\Rightarrow$(1),
  we proved that the image of $\CQK$ is connected and that $j$ is
  injective on the image of $\CQK$.

  Therefore $K_\CA$ is in the image of $\CQK$, and the map
  $j:T_{\CK_{\CA}}\to\bar T$ is injective. From the last paragraph of
  the proof of Theorem~\ref{thm:mainthm} we deduce that $j$ is
  isometric. Finally, from the minimality of $T$ we deduce that the
  image of $j$ contains $T$.
\end{proof}

\section{Applications to geometric trees and limits}
\label{sec:applications}

In this section we will present some first applications of the main
result of this paper, Theorem \ref{thm:mainthm}, to questions which in
part date back to the work of Gaboriau-Levitt \cite{gl}. It should
also be noted that Theorem \ref{thm:mainthm} is the basis for the
forthcoming papers \cite{coulbois-fractal} and \cite{ch}.

\smallskip

Recall that Outer space $\CVN$ is the space of projectivized minimal
free simplicial actions of $\FN$ on $\R$-trees. It comes with a
natural action by $\Out (\FN)$, and it is in many ways the analogue of
Teichm\"uller space, equipped with its action of the mapping class
group.  In particular, $\CVN$ has a natural ``Thurston boundary''
$\partial \CVN$, which defines a compactification $\barCVN = \CVN \cup
\partial \CVN$ of $\CVN$.  Its preimage $\barcvn$, obtained through
unprojectivization, consists precisely of all $\R$-trees $T$ with
non-trivial minimal very small action of $\FN$ by isometries.

\subsection{Geometric trees}

There is a special class of group actions on $\R$-trees which play an
important role in what is often called the ``Rips machine'': A minimal
$\R$-tree $T$ is called {\em geometric} if there exists a finite
subtree $K \subset T$ and a basis $\CA$ of $\FN$ such that the map $j:
T_{\CK} \to T$ is an isometry. It is proved in \cite{gl} that in this
case for any basis $\CA$ one can find such a finite subtree $K$.  For
more information about geometric trees regarding the context of this
paper see \cite{gl}.

Recall from \S\ref{limitset} that the heart $K_{\CA} \subset T$
denotes the convex hull of the limit set $\Omega_{\CA}$ of $T$ with
respect to the basis $\CA$.  We denote by $\CK_{\CA} = (K_{\CA}, \CA)$
the associated system of partial isometries.

\begin{cor}
\label{cor:geometrictree}
A very small minimal $\R$-tree $T$, with isometric $\FN$-action that
 has dense orbits, is geometric if and only if, for any basis $\CA$ of
 $\FN$, the heart $K_{\CA} $ is a finite subtree of $T$.
\end{cor}

\begin{proof}
If $T$ is geometric, then by definition there is a finite tree $K \subset 
T$ such that the map $j: T_{K} \to T$ is an isometry. Thus condition 
(1) of Theorem~\ref{thm:mainthm} is satisfied, and hence condition 
(3) implies that $K_{\CA}$ is a subtree of $K$, and thus it is finite.

Conversely, if $K_{\CA}$ lies in $T$, the image of the map $j$ defined
on $T_{\CK_{\CA}}$ is contained in $T \subset \bar T$, giving a map
$j: T_{\CK_{\CA}} \to T$ which by Corollary~\ref{cor:KA} is
isometric. By minimality of $T$, the map $j$ is onto.
\end{proof}

\subsection{Increasing systems of isometries}

Let $T$ be a minimal $\R$-tree with a very small action of $\FN$ by
isometries, which has dense orbits. As in \S\ref{sectionfive}, let
$\CA$ be a  basis of $\FN$, and for any $n \in \N$ let $K({n})$
be a compact subtree of the metric completion $\bar T$ of $T$, with non-empty intersections $K(n)
\cap a_{i} K(n)$ for all $a_{i}$ of $\CA$.  One obtains systems of
partial isometries $\CK(n) = (K(n), \CA)$ as in the previous sections.

We will consider sequences $K(n)$ which are increasing, i.e. for 
all $n \leq m$ we assume
\[
K(n) \subset K(m) \, .
\]
We can apply Theorem~\ref{thm:TK} to the case $K = K(n)$ and the tree
$T_{\CK(m)}$, to obtain canonical $\FN$-equivariant maps
\[
j_{m,n}:T_{\CK(n)} \to T_{\CK(m)}
\]
which satisfy $j_{k,m} \circ j_{m,n} = j_{k,n}$, for any natural
numbers $n \leq m \leq k$.

The maps $j_{m,n}$ are length decreasing morphisms, so that the trees
$T_{\CK(n)}$ converge in the equivariant Gromov-Hausdorff topology
(see \cite{paulin-top-equiv}) to an $\R$-tree $T_\infty$, equipped
with an action of $\FN$ by isometries.  Alternatively, one can pass to
the direct limit space defined by the system of maps $j_{m,n}$, which
inherits from the $T_{\CK(n)}$ a canonical pseudo-metric as well as an
action of $\FN$ by (pseudo-)isometries. One then defines $T_\infty$ as
the canonically associated metric quotient space.  Both,
arc-connectedness and 0-hyperbolicity carry over in those transitions,
so that $T_\infty$ is indeed an $\R$-tree with isometric $\FN$-action.

The minimal $\FN$-invariant subtrees $T_{\CK(n)}^{\mbox{\scriptsize
min}} \subset T_{\CK(n)}$ and $T_\infty^{\mbox{\scriptsize min}}
\subset T_\infty$ define points in the closure $\barcvn$ of
unprojectivized Outer space $\cvn$ (compare \cite{chl1-II} and the
references given there).  The sequence of trees
$T_{\CK(n)}^{\mbox{\scriptsize min}}$ converges in $\barcvn$ to the
tree $T_\infty^{\mbox{\scriptsize min}}$.

\smallskip

The maps $j_{m,n}$ also converge to $\FN$-equivariant maps $j_{\infty,
n}: T_{\CK(n)} \to T_{\infty}$ that satisfy $j_{\infty,m} \circ j_{m,
n} = j_{\infty,n}$.

\smallskip

We consider the increasing union of the $K(n)$, and we define
$K(\infty)$ to be its closure in $\bar T$,
\[
K(\infty) = \overline{\bigcup_{n \in \N} K(n)}\, ,
\]
provided with the induced system $\CK(\infty) = (K(\infty), \CA)$ of
partial isometries.  We always assume that $K(\infty)$ is compact.

\smallskip

Using that $K(n)\subset  K(\infty)$, we can apply again
Theorem~\ref{thm:TK} to get $\FN$-equivariant, length decreasing
morphisms:
\[
j_{0,n}:T_{\CK(n)} \to T_{\CK(\infty)}
\]
These maps converge to an $\FN$-equivariant length decreasing map
$j_{0,\infty}:T_\infty \to T_{\CK(\infty)}$. This map continuously
extends to the metric completions $\bar{\jmath}_{0,\infty}:\bar
T_\infty \to \bar T_{\CK(\infty)}$.

\smallskip

For each $n \in \N$, the map $j_{0,n}$ restricts to an isometry on
$K(n)$, and thus $j_{0,\infty}$ restricts to an isometry on the union
of the $K(n)$, which extends to an isometry from $K(\infty)\subset
\bar T_\infty$ onto its image in $T_{\CK(\infty)}$.  Applying
Theorem~\ref{thm:TK} again, to the inverse of this isometry, we get an
$\FN$-equivariant length decreasing morphism $j_{\infty,0}:
T_{\CK(\infty)} \to \bar T_\infty$.

\smallskip

By construction, the restrictions of the maps $j_{0,\infty}$ and
$j_{\infty,0}$ to each of the $K(n)$ are isometries which are inverses
of one another. Thus the map $\bar{\jmath}_{0,\infty}\circ
j_{\infty,0}$ is length decreasing and restricts to the identity on
$\underset{n\in \N}{\bigcup} K(n)$ and thus on $K(\infty)$. Using
Theorem~\ref{thm:TK}, we see that it is an isometry on all of
$T_{\CK(\infty)}$.  This shows
\[
\lim_{n\to +\infty} \, T_{\CK(n)} = T_\infty \subset T_{\CK(\infty)}
\]
and thus
\[
\lim_{n\to +\infty} \, T^{\min}_{\CK(n)} = T^{\min}_{\CK(\infty)} \, .
\]

\smallskip

As a direct consequence of Theorem~\ref{thm:mainthm} one now derives:

\begin{cor}
\label{cor:approxM}
Let $T$ be a minimal $\R$-tree with a very small action of $\FN$ by
isometries, which has dense orbits. Let $\CA$ be a basis of $\FN$. For
any $n \in \N$, let $K({n})$ be a compact subtree of $\bar T$ with
non-empty intersections $K(n) \cap a_i K(n)$, for all $a_i$ of $\CA$.
Let $\CK(n) = (K(n), \CA)$ be the induced systems of isometries. Let
$K(\infty)$ be the closure of the increasing union of the $K(n)$, and
assume that $K(\infty)$ is compact.

Then the minimal trees $T_{\CK(n)}^{\mbox{\scriptsize min}}$ converge
in $\barcvn$ to $T$ if and only if $K(\infty)$ contains
$\Omega_{\CA}$.\qed
\end{cor}

An application of this corollary is the following sharpening of a
classical result of Gaboriau-Levitt \cite{gl}, who showed that every
$T \in \barcvn$ can be approximated by a sequence of geometric
$T_{\CK(n)}$, i.e. each $K(n)$ is a finite subtree of $T$.

\begin{cor}
\label{cor:finite-bounded-approximation}
For every very small minimal $\R$-tree $T$, with isometric
$\FN$-action that has dense orbits, there exists a sequence of finite
subtrees $K(n)$ of uniformely bounded diameter, such that:
\[
T = \lim_{n \to \infty} \, T_{\CK(n)}
\]
\end{cor}

\begin{proof}
It is well known \cite{gl} that the number of branch points in $T$ is
a countable set $\cal P$ that is dense in every segment of $T$.  It
suffices to consider the countable subfamily $(P_n)_{n\in \N} = {\cal
P} \cap K_\CA$ and to define $K(n)$ as convex hull of the set $\{P_1,
\ldots, P_n\}$. Since the heart $K_\CA$ is the convex hull of the
limit set $\Omega_\CA$, the claim is a direct consequence of Corollary
\ref{cor:approxM}.
\end{proof}

\subsection{Approximations by simplicial trees}

An algebraic lamination $L$ is said to be {\em closed by diagonal
leaves}, if for any leaves $(X,X')$ and $(X',X'')$ in $L$ one either
has $X=X''$, or $(X,X'')$ is again a leaf in $L$.  We remark that, if
$T$ is an $\R$-tree with a minimal action of $\FN$ by isometries that
has dense orbits, it follows from \S\ref{themapQ2} that the dual
lamination $L(T)$ of $T$ is closed by diagonal leaves. Also, for a
system of isometries $\CK$ with independent generators, we deduce from
\S\ref{admissiblelamination} and \S\ref{independent-generators} that
the admissible lamination $\Ladm(\CK)$ is closed by diagonal leaves.

An algebraic lamination $L$ is said to be {\em minimal up to diagonal
leaves} if it does not contain a proper non-trivial sublamination
that is closed by diagonal leaves.

\begin{cor}
\label{cor:approximationfree}
Let $T$ be an $\R$-tree with a minimal, very small action of $\FN$
that has dense orbits. Let $\CA$ be a basis of $\FN$. If a compact
subtree $K \subset \bar T$ does not contain $\Omega_{\CA}$, and if the
dual lamination $L(T)$ is minimal up to diagonal leaves, then the
approximation tree $T_{\CK}^{\mbox{\scriptsize min}}$ is free simplicial (i.e. it belongs to
the unprojectivized Outer space $\cvn$ rather than to its boundary $\partial\cvn$).
\end{cor}

\begin{proof}
  From Proposition~\ref{prop:dualvsadm} we know that the dual
  lamination $L(T_\CK)$ of $T_\CK$ is a sublamination of the
  admissible lamination $\Ladm(\CK)$.  The lamination $\Ladm(\CK)$ is
  closed by diagonal leaves and is a sublamination of $L(T)$, by
  Lemma~\ref{lem:lamsubset}. Since $\Omega_\CA$ is not a subset of $K$,
  Theorem~\ref{thm:mainthm} implies that the admissible lamination $\Ladm(\CK)$ is a strict sublamination of
  $L(T)$.

From the minimality of $L(T)$ up to diagonal leaves we deduce that
$\Ladm(\CK)$ and $L(T_\CK)$ are empty, so that (compare
\S\ref{duallamination}) the action of $\FN$ on $T^{\mbox{\scriptsize
min}}_\CK$ is free and discrete.
\end{proof}

This corollary indicates that
the resolution of an arbitrary $\R$-tree with
isometric $G$-action, for more general groups $G$, via systems of
partial isometries on a finite tree, as promoted by the Rips machine, may yield directly
a simplicial tree,
i.e. without having to go through further iterations in Rips' procedure.

\bigskip

\noindent {\em Acknowledgments:} The authors would like to thank
V.~Guirardel, P.~Hubert and G.~Levitt for helpful comments.  The first
and the third author would also like to thank the MSRI at Berkeley for
the support received from the program ``Geometric Group Theory'' in the
fall of 2007.

%\bibliographystyle{plain}
%\bibliography{../CHL/bibli.bib}
%\end{document}

\end{document}